# RESULTS ON FIXED POINTS OF CLOSED AND COMPACT COMPOSITE SEQUENCES OF OPERATORS AND PROJECTORS IN A CLASS OF COMPLETE METRIC SPACES

M. De la Sen, IIDP. Faculty of Science and Technology. University of the Basque Country

Campus of Leioa (Bizkaia). PO Box 644- Bilbao, SPAIN

**Abstract.** Some results on fixed points related to the contractive compositions of bounded operators in complete metric spaces are discussed through the manuscript. The class of composite operators under study can include, in particular, sequences of projection operators under, in general, oblique projective operators.

## 1. Introduction

Some results on fixed points related to the contractive compositions of bounded operators in a class of complete metric spaces $(X, d)$ which are also Banach spaces if $X$ is a vector space on a certain field $F$ (usually $R$ or $C$) and the metric is homogeneous and translation-invariant, so that it is also a norm, are discussed through the manuscript. The class of composite operators under study can include, in particular, sequences of projection operators under, in general, oblique projective operators. Section 2 is concerned of composite operators which include sequences of pairs of contractive operators including, in general, oblique projection operators in the operator composite strip. The results are generalized in Section 3 to sequences of, in general, non constant bounded closed operators which can have bounded, closed and compact limits and such that the relevant composite sequences are also compact operators. This manuscript addresses Banach´s contraction principle [1-4] guaranteeing the existence of unique fixed points under contractive conditions fulfilled by some relevant strips of composite operators within in the whole composite sequence of operators.

## 2. Some results on contractive mappings and fixed points under projection operators

Let $\{T_k\}$ be a sequence $T_k : X \to X$ of self-mappings on a vector space $X$ over a field $F$, where $(X, d)$ is a metric space and consider a sequence $\{P_k\}$ of (non-necessarily orthogonal) projection operators on $X$ of respective ranges $M_k$ which are then closed subspaces of $X$, [3]. We can then consider a sequence of projection operators $\{P_{M_k}\}$ with $P_{M_k} : X \to M_k$ such that $P_k = P_{M_k}$ so that $X = Im P_k \oplus Ker P_k$ and $z = P_k x \in Im P_k$ is in $M_k$ for any $x \in X$ and $\bar{z} = x - z = (I - P_k)x \in Ker P_k$ for $k \in N_0 = N \cup \{0\}$, Now, consider sequences $\{x_k\}$ in $X$ and $\{z_k\}$ in $M_k$ with $z_k = P_k x_k$ such that the identities:

$$x_{k+1} = T_k x_k = P_{k+1} x_{k+1} + (I - P_{k+1}) x_{k+1} = z_{k+1} + (I - P_{k+1}) x_{k+1} = T_k P_k x_k + T_k (I - P_k) x_k \quad (2.1)$$

hold by construction for $k \in N_0$. The subsequent result holds:

**Theorem 2.1**. Assume that $(X, d)$ is a complete metric space with the metric $d: X \times X \to \mathbf{R}_{0+}$ being homogeneous and translation- invariant and $0 \in X$ ; $k \in \mathbf{N}_0$. The following properties hold:

**(i)** If all the self- mappings on $X$ in the sequence $\{T_k\}$ are non-expansive and the sequence of projection operators $\{P_k\}$ from $X$ to the sequence of subspaces $\{M_k\}$ is uniformly bounded then $d(z_{k+2}, z_{k+1}) < \infty$ ; $k \in \mathbf{N}_0$.

**(ii)** Assume that the self- mappings on $X$ in the subsequence $\{T_k\}_{k \geq n_0}$ are contractive for some $n_0 \in \mathbf{N}_0$, that the sequence of operators $\{T_k\}$ converges to $T: X \to X$ and that the projection operator $P: X \to M \subset X$ is constant and bounded (i.e. if it is not orthogonal, that is, it is oblique, then its norm exceeds one and it is finite) then Property (i) holds. Furthermore, $\exists \lim_{k \to \infty} d(z_{k+2}, z_{k+1}) = 0$ and $\{z_{k+1} = PT_k x_k\}$ is a Cauchy sequence which converges to some unique limit point $z(= PTx) \in M$ for any initial iterate $x_0 \in X$ where $x(= Tx) \in X$ is the unique fixed point of $T: X \to X$.

**(iii)** Assume that there is a strictly sequence of nonnegative integers $\{j_k\}$ such that the difference sequence $\{j_{k+1} - j_k\}$ is uniformly bounded and has a limit $J \in \mathbf{N}$ as $k \to \infty$. Assume also that the associate sequence of composite self-mappings $\hat{T}(j_{k+1} + 1, j_k) = \{T_{j_{k+1}} \ldots T_{j_k+1} T_{j_k}\}$ is contractive and that the sequence projection operators $\{P_k\}$ from $X$ to the sequence of subspaces $\{M_k\}$ is uniformly bounded and has a set of subsequences each converging to a set $\{P_{j_i}\} \to \hat{P}_i$ of projectors from $X$ to $\{M_{j_i}\}$. Then, $\exists \lim_{k \to \infty} d(z_{j_{k+1}+1}, z_{j_k+1}) = 0$, and there is at most a finite number $J$ of distinct Cauchy subsequences $\{z_{j_k}\}$ with distinct limit points $\{\hat{z}_1, \ldots, \hat{z}_J\}$ in $X$.

*Proof*: Since the metric is homogeneous and translation -invariant then the complete metric space $(X, d)$ associated with the vector space $X$ can also be considered as a Banach space $(X, \|\ \|)$ under the metric-induced norm defined as $\|x - y\| = d(x, y)$ ; $\forall x, y \in X$. The norm of any projection operator in the considered sequence is defined as $\|P_k\| = \sup_{\|x(\in X) \neq 0\| \leq 1} \frac{\|P_k x\|}{\|x\|} = \sup_{\|x(\in X)\| = 1} \frac{\|P_k x\|}{\|x\|} = \sup_{\|x(\in X)\| = 1} \|P_k x\|$ ; $k \in \mathbf{N}_0$. Then, if $\hat{T}(k+1, j): X \to X$ is the left-composite self-mapping $\hat{T}(k+1, j) = T_k \ldots T_{j+1} T_j$ for any $k(\geq j)$, $j \in \mathbf{N}_0$, one gets from direct calculations, by using the property that the metric is homogeneous and translation-invariant, the following relations for any iterated sequences $\{x_k\}$ in $X$ and $\{z_k\}$ in $M_k$ constructed as $x_{k+1} = T_k x_k$, $z_k = P_k x_k$ with arbitrary $x_0 \in X$ ; $k \in \mathbf{N}$ :

$$d(z_{k+2}, z_{k+1}) = d(P_{k+2} x_{k+2}, P_{k+1} x_{k+1}) = d(P_{k+1} x_{k+2}, P_{k+1} x_{k+1} + P_{k+1} x_{k+2} - P_{k+2} x_{k+2})$$
$$\leq d(P_{k+1} x_{k+2}, P_{k+1} x_{k+1}) + d(P_{k+1} x_{k+1}, P_{k+1} x_{k+1} + P_{k+1} x_{k+2} - P_{k+2} x_{k+2})$$

$$\leq d(P_{k+1}x_{k+2}, P_{k+1}x_{k+1}) + d(P_{k+2}x_{k+2}, P_{k+1}x_{k+2})$$

$$= \|P_{k+1}x_{k+2} - P_{k+1}x_{k+1}\| + \|P_{k+2}x_{k+2} - P_{k+1}x_{k+2}\|$$

$$\leq \|P_{k+1}\| \|x_{k+2} - x_{k+1}\| + \|P_{k+2} - P_{k+1}\| \|x_{k+2}\|$$

$$= \|P_{k+1}\| d(x_{k+1}, x_{k+2}) + \|P_{k+2} - P_{k+1}\| d(x_{k+2}, 0)$$

$$= \|P_{k+1}\| d(T_k x_k, T_{k+1} T_k x_k) + \|P_{k+2} - P_{k+1}\| d(\hat{T}(k+2, 0) x_0, 0)$$

$$\leq \|P_{k+1}\| d(x_1, x_0) + \|P_{k+2} - P_{k+1}\| d(\hat{T}(k+2, 0) x_0, 0)$$

$$\leq \|P_{k+1}\| d(x_1, x_0) + (\|P_{k+1}\| + \|P_{k+2}\|) d(x_0, 0)$$

$$\leq \|P_{k+1}\| (d(T_0 x_0, 0) + d(x_0, 0)) + (\|P_{k+1}\| + \|P_{k+2}\|) \|x_0\| \; ; \; k \in N_0 \quad (2.2)$$

since the metric is homogeneous and translation-invariant, the norm is an induced-metric norm, then $\|x_k\| = d(x_k, 0) = d(\hat{T}(k, 0)x_0, 0)$ where $z_k \in M_k$, and the self-mappings in the sequence $\{T_k\}$ on $X$ are all non-expansive; $k \in N_0$, and the sequence of projection operators $\{P_k\}$ from $X$ to the sequence of subspaces $\{M_k\}$ is uniformly bounded with $\sup_{k \in N_0} \|P_k\| \leq \mu < \infty$. Then, one has from (2.2):

$$d(z_{k+2}, z_{k+1}) \leq 4\mu \|x_0\| < \infty \; ; \; k \in N_0 \quad (2.3)$$

where $\mu = 1$ if all the projections are orthogonal and $\mu > 1$, otherwise. Hence, Property (i). If $P_k = P_{k+1} = P$ for $k \in N_0$ is a constant bounded projection from $X$ to $M$ with $M_k = M \subset X$ being constant for $k \in N_0$ and all the self-mappings on $X$ of the sequence $\{T_k\}$ are contractive then one gets from (2.2) for the real constant $K = \sup_{k(\geq n_0) \in N_0} K_k$ such that $K \in [0, 1)$ that Property (i) holds according to the relation:

$$d(z_{k+2}, z_{k+1}) \leq \|P\| d(T_k x_k, T_{k+1} x_{k+1}) \leq \|P\| d(T_k x_k, T_k x_{k+1} + T_{k+1} x_{k+1} - T_k x_{k+1})$$

$$\leq \|P\| d(T_k x_k, T_k x_{k+1}) + \|P\| d(T_k x_{k+1}, T_k x_{k+1} + T_{k+1} x_{k+1} - T_k x_{k+1})$$

$$= \|P\| (d(T_k x_k, T_k x_{k+1}) + \|P\| d(0, T_{k+1} x_{k+1} - T_k x_{k+1}))$$

$$= \|P\| (d(T_k x_k, T_k x_{k+1}) + d(T_k x_{k+1}, T_{k+1} x_{k+1}))$$

$$= \|P\| (\|T_k x_k - T_k x_{k+1}\| + \|T_k x_{k+1} - T_{k+1} x_{k+1}\|)$$

$$\leq \|P\| \|T_k\| \|x_k - x_{k+1}\| + \|P\| \|T_k - T_{k+1}\| \|x_{k+1}\|$$

$$\leq \|P\| (\|T_k\| \|x_k - x_{k+1}\| + \|T_k - T_{k+1}\| d(\hat{T}(k+1, 0) x_{n_0}, \hat{T}(k+1, 0) 0))$$

$$\leq \|P\| \left( K^{k-n_0+1} d(x_{n_0+1}, x_{n_0}) + \frac{K(1 - K^{k-n_0})}{1 - K} \|T_k - T_{k+1}\| \|x_{n_0}\| \right) < \infty \; ; \; k(\geq n_0) \in N_0 \quad (2.4)$$

so that $\exists \lim\limits_{k \to \infty} d(z_{k+2}, z_{k+1}) = 0$ from (2.4) for any initial value $x_0 \in X$ of the iteration since $K^{k-n_0-1} \to 0$ and $\|T_k - T_{k+1}\| \to 0$ as $k \to \infty$ since $\{T_k\} \to T$. Then, $\{z_k\}$ is Cauchy sequence which has a limit $z$ in $M$, since $M$ is closed, [4]. It is now proven that $z(\in M) = Px$ is the unique limit point in $M$ of any sequence of iterates where $x(=Tx) \in X$ is a fixed point of the self-mapping $T: X \to X$ which is unique from Banach contraction principle. It is now proven that $T: X \to X$ is contractive. Assume not so that one has if $T: X \to X$ is not contractive:

$$d(x_{k+1}, x_k) \leq d(Tx_{k+1}, Tx_k) = d(T_k x_{k+1}, T_k x_k + Tx_k - T_k x_k + T_k x_{k+1} - Tx_{k+1})$$
$$\leq d(T_k x_{k+1}, T_k x_k) + d(T_k x_k, T_k x_k + Tx_k - T_k x_k + T_k x_{k+1} - Tx_{k+1})$$
$$= d(T_k x_{k+1}, T_k x_k) + d(0, Tx_k - T_k x_k + T_k x_{k+1} - Tx_{k+1})$$
$$= d(T_k x_{k+1}, T_k x_k) + d(T_k x_k + Tx_{k+1}, Tx_k + T_k x_{k+1})$$
$$\leq d(T_k x_{k+1}, T_k x_k) + \|T_k - T\|(\|x_k\| + \|x_{k+1}\|)$$
$$\leq K_k d(x_{k+1}, x_k) + \|T_k - T\|(\|x_k\| + \|x_{k+1}\|); \; k \geq n_0$$

for nonzero $x_k$ and $x_{k+1}$ since $T_k: X \to X$ is contractive for $k \geq n_0$. Then, on gets, since $\|T_k - T\| \to 0$ as $k \to \infty$, that

$$\lim\sup_{k \to \infty} \left[ (1 - K_k - \|T_k - T\|)(\|x_k\| + \|x_{k+1}\|) \right] = \lim\sup_{k \to \infty} \left[ (1 - K_k)(\|x_k\| + \|x_{k+1}\|) \right] \leq 0$$

which is a contradiction since $K_k < 1$ for $k \geq n_0$ unless $\{x_k\}$ converges to zero. If $\{x_k\}$ converges to zero then there are $n_1(\geq n_0) \in N_0$ and $0 < \lambda = \lambda(n_1) < 1 - \sup\limits_{k \geq n_1} K_k$ such that $\|T_k - T\| \leq \lambda$ for all $k \geq n_1$ since $\|T_k - T\| \to 0$ as $k \to \infty$ and some $k_1 \geq n_1$ such that $\|x_{k_1}\| + \|x_{k_1+1}\| > 0$ that yields the contradiction

$$0 < (1 - K_{k_1} - \lambda)(\|x_{k_1}\| + \|x_{k_1+1}\|) \leq 0$$

Thus, if the subsequence $\{T_k\}_{k \geq n_0}$ is contractive in $X$ then its limit $T: X \to X$ is also contractive. Now, since $T: X \to X$ is contractive then its fixed point is unique since $(X, d)$ is complete. It is clear that $z = PTx$ is a limit point in $M$ of any iterated sequence. Assume that it is not unique so that there are two limit points $z = PTx$, $\hat{z} = PT\hat{x}(\neq z) \in M$ for some $\hat{x}(\neq x) \in X$ which is not trivially a fixed point of $T: X \to X$ (since the fixed point $x \in X$ of the contractive self-mapping $T: X \to X$ is unique if $(X, d)$ is complete). Thus, from Banach contraction principle and since $(X, d)$ is complete, one has

$$0 \leftarrow d(PT_k \hat{x}, PT_k x) = \|PT_k \hat{x} - PT_k x\| \leq \|P\| \|T_k \hat{x} - T_k x\| \leq \|P\| K^k d(\hat{x}, x) < \infty \; ; \; k \in N_0 \qquad (2.5)$$

as $k \to \infty$ since $\{T_k\}$ converges, there is a limit self-mapping $T$ on $X$:

$$\hat{z} \leftarrow PT_k \hat{x} \to PT_k x = PTx = Px = z \qquad (2.6)$$

Thus, $PT_k x \to \hat{z} = z$ as $k \to \infty$. Hence a contradiction to $\hat{z} \neq z$ and then $z$ in $M$ is the unique limit point of $PT: X \to M$ even in the event that there is $\hat{x}(\neq x = Tx) \in X$ such that $PT\hat{x} = PTx = Px = z$. Property (ii) has been proven.

On the other hand, if the sequence of operators is uniformly bounded then $\sup_{k \in N_0} \|P_k\| \leq \mu < \infty$ and, if furthermore, the sequence of compositie mappings $\{\hat{T}(j_{k+1}, j_k)\}$ is contractive with some constant $\hat{K} \in [0,1)$ given by $\hat{K} = \sup_{k \in N_0} \left(\prod_{j=j_k}^{j_{k+1}-1}[K_j]\right)$, $k \in N_0$, where $\{j_k\}$ is a strictly increasing sequence of natural numbers such that the sequence $\{j_{k+1} - j_k\}$ is uniformly bounded, one has directly from (2.1)-(2.2):

$$d(z_{j_{k+1}}, z_{j_k}) \leq \mu\, d(x_{j_k}, \hat{T}(j_{k+1}, j_k)x_{j_k})$$
$$\leq \mu\left(\hat{K}^{j_k}[d(x_1, x_0) + 2\|x_0\|]\right) < \infty\; ; \; k \in N_0 \quad (2.7)$$

$$d(z_{j_k + i_k}, z_{j_k}) \leq \mu\, d(x_{j_k}, \hat{T}(j_{k+1}, j_k)x_{j_k})$$
$$\leq \mu\left(\prod_{j=j_k}^{j_k+i_k-1} K_j\right)\left(\hat{K}^{j_k}[d(x_1,x_0)+2\|x_0\|]\right) < \infty \quad (2.8)$$

for $k \in N_0$, where $\{\{i_k\}\}$ is a sequence of finite sets of natural numbers satisfying $j \in \{i_k\} \Leftrightarrow 1 \leq j < (j_{k+1} - j_k) + 1$ for $k \in N_0$. Thus, one gets from (2-7)-(2.8):

$$\exists \lim_{j_k \to \infty} d(z_{j_{k+1}}, z_{j_k}) = \lim_{j_k \to \infty} d(z_{j_{k+1}+j}, z_{j_k}) = 0 \quad (2.9)$$

for any $j \in \{i_k\}$ since $j_k \to \infty$ as $k \to \infty$ and. Since $\{j_{k+1} - j_k\}$ is uniformly bounded with existing limit $J \in N$ as $k \to \infty$, $\hat{T}(j_{k+1}+1, j_k) = \{T_{j_{k+1}} \ldots T_{j_k+1} T_{j_k}\}$ is contractive and the projection sequence $\{P_k\}$ from $X$ to the sequence of subspaces $\{M_k\}$ is uniformly bounded while having a finite set of subsequences $\{P_{j_i}\}$ converging to a set $P_{j_i} \to \hat{P}_i$ of projectors from $X$ to $\{M_k\}$ for $i \in \bar{J}$ then

$$\exists \lim_{k \to \infty} d(z_{j_{k+1}+1}, z_{j_k+1}) = \lim_{k \to \infty} d\left(P_{j_{k+1}+1}\hat{T}(j_{k+1}+1, 0)x_0, P_{j_k}\hat{T}(j_k+1, 0)x_0\right) = 0 \quad (2.10)$$

so that $\{z_{j_{k+1}+1} = P_{j_{k+1}}\hat{T}(j_{k+1}+1, 0)x_0\}$ is a Cauchy sequence satisfying:

$$z_{j_k+i} \to \hat{z}_i = \hat{P}_i\left(\lim_{k \to \infty}(T_{j_k+i} \ldots T_{j_{k+i}-1}T_{j_k})\hat{z}_{i-1}\right) \text{ as } j_i \to \infty\; i \in \bar{M}\;;\; \hat{z}_i \in \left(X \cap \bigcap_{k=0}^{\infty} M_{j_k+i-1}\right)$$

(2.11)

for at most $J$ distinct points $\{\hat{z}_1, \ldots, \hat{z}_J\}$ since, by hypothesis, there is a natural number $J$ satisfying $\infty > J = \lim_{k \in N_0} \{j_{k+1} - j_k\} \geq 1$ and since $\{P_{j_i}\} \to \hat{P}_i$; $i \in \bar{J}$. Hence, Property (iii) holds. □

*Remark 2.2.* The existence of some $\hat{x}(\neq x = Tx) \in X$ in the proof of Theorem 2.1(ii) often happens. For instance if $PT: X \to M$ is linear then $\hat{x} = x + x_a$ fulfils the relations $PTx = Px = PT\hat{x} = z$ for any $x_a \in \text{Ker } PT$. □

The following auxiliary result to be then used holds:

**Lemma 2.3.** Assume that the sequences of linear self-mappings $\{T_n\}$ and $\{P_n\}$ converge to respective limits $P$ and $T$ being mappings from $X$ to $M \subset X$ and from $X$ to $X$, respectively in the sense that $P_n T_n x \to PTx$ as $n \to \infty$ from any $x \in X$. Let $(X, \|\cdot\|)$ be a Banach space with the norm of any $Q: X \to X$ being defined by

$$\|Q\| = \sup_{x \in X, \|x\|=1} \|Qx\| \qquad (2.12)$$

For any given $\delta \in \mathbf{R}_+$, $\exists n_0 = n_0(\delta)$ such that $\|P_n T_n\| \leq pt + \delta$ where $p = \|P\|$ and $t = \|T\|$. If $PT: X \to X$ is contractive then the sequence $\{P_n T_n\}$ of mappings from $X$ to $X$ is then asymptotically contractive. □

## 3. Results on contractive mappings of sequences of composite bounded operators

The results of Theorem 2.1 are now extended to the study of contractive compositions of linear operators belonging to two sequences of bounded operators $\{T_{ki}\}$ with $T_{ik}: Dom(T_{ik}) \subset X \to Im(T_{ik}) \subset X$ to $X$; $i = 1,2$ so that none of them is necessarily a projection on some subspace of $X$. Some preparatory results are first established. In the following, a Banach space $(X, \|\cdot\|)$, where $X$ is a vector space on a field $F$, being equivalent to a complete metric space $(X, d)$ with a homogeneous and translation-invariant metric induced-norm $d: X \times X \to \mathbf{R}_{0+}$ is considered such that $\|x\| = d(x, 0) = d(x + \alpha y, \alpha y)$ for any real $\alpha$ and any $x, y \in X$. The subsequent result refers to the asymptotic distances in sequences involving a convergent composite sequence of bounded linear operators.

**Lemma 3.1.** Consider a Banach space $(X, \|\cdot\|)$, with $0 \in X$, being equivalent to a complete metric space $(X, d)$ with a homogeneous and translation-invariant metric induced-norm $d: X \times X \to \mathbf{R}_{0+}$. Consider also a composite sequence of two sequences of bounded linear operators $\{T_k = T_{2k} T_{1k}\}$ defined by $T_k: Dom(T_k) \subset X \to Im(T_k) \subset X$ defined by $T_k x = T_{2k}(T_{1k} x) = T_{2k} T_{1k} x$ for any $x \in Dom(T_k)$, where $T_{ik}: Dom(T_{ik}) \subset X \to Im(T_{ik}) \subset X$; $i = 1,2$, provided that $Im(T_{1k})(\neq \emptyset) \subset Dom(T_{2k})$ and $Im(T_{2k}) \cap Dom(T_{1k+1}) \neq \emptyset$; $\forall k \in \mathbf{N}_0$. The following properties hold:

**(i)** Assume that $\{T_{ik}\} \to T_i$ $(i = 1,2)$. Then, $\lim_{k \to \infty} d(T_k x, T_2 T_1 x) = 0$; $\forall x \in Dom(T_k)$; $\forall k \in \mathbf{N}_0$.

**(ii)** Assume that $\{T_{1k}\} \to T_1$ $(i = 1,2)$. Then, $\lim_{k \to \infty} d(T_k x, T_{2k} T_1 x) = 0$; $\forall x \in Dom(T_k)$; $\forall k \in \mathbf{N}_0$.

**(iii)** Assume that $\{T_{2k}\} \to T_2$ $(i=1,2)$. Then, $\lim_{k \to \infty} d(T_k x, T_2 T_{1k} x) = 0$; $\forall x \in Dom(T_k)$; $\forall k \in N_0$.

**(iv)** Define the operator composite sequence $\{\hat{T}(k+i+1,k)\}$ of operators as $\hat{T}(k+i+1,k) = T_{k+i}...T_{k+1}T_k$ with $T_k = T_{j_k k} T_{j_k-1,k}....T_{2k} T_{1k}$; $\forall i, k \in N_0$ subject to $Im(T_{jk})(\neq \varnothing) \subset Dom(T_{j+1,k})$ for $j \in \overline{j_k}$, $ImT_{j_k k}(\neq \varnothing) \subset DomT_{1,k+1}$ and $1 \leq j_k \leq J(\in N_0) < \infty$; $\forall k \in N_0$. Define also the operator composite sequence of operators $\{\hat{T}^0(k+i+1,k)\}$ as

$$\hat{T}^0(k+i+1,k) = T^0_{k+i}...T^0_{k+1}T^0_k = \left(T^0_{j_{k+i} k+i}...T^0_{j_{k+i}-1,k+i} T^0_{1,k+i}\right)....\left(T^0_{j_k k}...T^0_{2k} T^0_{1k}\right) \quad (3.1)$$

for $i, k \in N_0$, where $T^0_{jk} = T_{jk}$ if $T_{jk}$ has not a limit as $k \to \infty$ and $T^0_{jk} = T_{j\infty}$ if $T_{j\infty} = \lim_{k \to \infty} T_{jk}$. Then

$$\lim_{k \to \infty} d\left(\hat{T}(k+i+1,i)x, \hat{T}^0(k+i+1,i)x\right) = 0 \; ; \; \forall x \in Dom(T_{10}), \; \forall i \in N_0$$

Properties (i)-(iii) hold for any $x \in Dom(T_{10})$. □

Property (iv) is direct from Properties (i)-(iii) and the associative property of composition of operators since for any $k \in N_{0+}$, $T_k x$ exists in $Im(T_k)$ if $Im(T_{1k}) \subset Dom(T_{2k})$ and $Dom(T_{1k+1}) \supset Im(T_{2k})$ for $x \in Dom(T_{10})$, and then, $Dom(T_k) \subset Dom(T_{1k})$ and $Im(T_k) \supset Im(T_{2k})$; $\forall k \in N_0$ since

$$\hat{T}(k+i+1,k) = T_{k+i}...T_{k+1}T_k = \left(T_{j_{k+i} k+i}...T_{j_{k+i}-1,k+i} T_{1,k+i}\right)....\left(T_{j_k k}...T_{2k} T_{1k}\right); \; \forall k \in N_0 \quad (3.2)$$

Then, for any finite $i \in N$, one gets:

$$\left\|\hat{T}(k+i+1,k) - \hat{T}^0(k+i+1,k)\right\| \leq \sum_{j=k}^{k+i} c_j \left[o\left(\left\|T^0_j - T_j\right\|\right)\right]$$

$$\leq \sum_{j=k}^{k+i} \sum_{\ell_j=1}^{j_k} c_{j\ell_j} \left[o\left(\left\|T^0_{\ell_j j} - T_{\ell_j j}\right\|\right)\right] \to 0 \text{ as } k \to \infty \quad (3.3)$$

for some positive finite constants $c_j$ and $c_{j\ell_j}$ since any linear operator $T^0_{\ell_j j}$ with a limit $T_{\ell_\infty j}$ admits a unique decomposition $T^0_{\ell_j j} = T_{\ell_\infty j} + \tilde{T}_{l_j j}$, with $\tilde{T}_{l_j j} \to 0$ as $k \to \infty$, $\forall \ell_j \in \overline{j_k}$, $j = 1, 2, ..., k+i$.

The subsequent result given without proof is concerned with the closeness of the limit operator if the sequence of operators is closed.

**Lemma 3.2**. Consider a sequence of closed linear operators $\{T_n\}$ defined by $T_n : Dom(T_n) \subset X \to Im(T_n) \subset X$ in a Banach space $(X, \|\;\|)$, such that $Im(T_n) \subset Dom(T_{n+1})$ with $Im(T_n) \cap Dom(T_{n+1}) \neq \varnothing$, which converge to a limit operator $T: Dom(T) \subset X \to Im(T) \subset X$. Then, such a limit is a closed operator which is bounded if all the operators of the sequence are bounded. □

The limit of bounded converging sequences belongs to the domain of the limit operator. Furthermore, one has for any bounded sequence $\{x_n\}$ converging to $x \in Dom(T)$:

$$\|T_n x_n - Tx\| = \|Tx_n + (T_n - T)x_n - Tx\| \le \|T\|\|x_n - x\| + \|T_n - T\|\|x_n\| \to 0 \qquad (3.4)$$

as $k \to \infty$ then $T_n x_n \to Tx$ strongly so that $T: Dom(T) \subset X \to Im(T) \subset X$ is a closed operator as a result.

The above result can be extended to sequences of operators not all of them being bounded provided that each of such sequences of operators can be decomposed as a composition of subsequences of composite operators such that each of such a composite subsequence is bounded. The above result can be applied to sequences of operators not all of them being bounded. It is well-known that a sequence of linear operators on a Hilbert space [5-6] is bounded if and only if they are closed and their domain is the whole vector space $X$, [1, 4]. Thus, we can obtain following result from Lemma 3.2:

**Lemma 3.3**. Consider a sequence of linear bounded operators $\{T_n\}$ defined by $T_n: X \to X$ in a Banach space $(X, \|\ \|)$ which converge to a limit operator $T: X \to X$. Then, such a limit is a bounded linear (and then continuous and closed) operator. □

The subsequent result, stated without proof, is concerned with the limit operator of a sequence of linear operators being compact if all the operators in the sequence are bounded and at least one of them is compact:

**Lemma 3.4**. The following properties hold:

**(i)** Consider a sequence of bounded compact linear operators $\{T_n\}$ defined by $T_n: Dom(T_k) \subset X \to Im(T_n) \subset X$ in a Banach space $(X, \|\ \|)$, such that $Im(T_n) \subset Dom(T_{n+1})$ with $Im(T_n) \cap Dom(T_{n+1}) \ne \emptyset$, which converge to a limit operator $T: Dom(T) \subset X \to Im(T) \subset X$. Then, such a limit is a compact operator.

**(ii)** Assume that the sequence $\{T_n\}$ of bounded operators satisfies that there is at least one compact operator within all subsequences $\{T_{j_n}, T_{j_n+1}, \ldots, T_{j_{n+1}-1}\}$ being subject to $\max_{n \in N_0}(j_{n+1} - j_n) \le c_j < \infty$ for some subsequence $\{j_n\} \subset N_0$ for any $n \in N_0$. Then, the composed operator $\hat{T}(n,m)$ is compact as it is its limit provided that it exists.

*Proof*: We have to prove that if $\{x_n\}$ is bounded then $\{Tx_n\}$ is convergent. Note that for given bounded sequences $\{x_n^{(i)}\}$ and $\{x_n^{(j)}\}$; $i, j, n \in N_0$ that

$$\|T x_n^{(i)} - T x_n^{(j)}\| = \|(T - T_n)x_n^{(i)} + T_n x_n^{(i)} - T_n x_n^{(j)} - (T - T_n)x_n^{(j)}\|$$

$$\le \|(T - T_n)x_n^{(i)}\| + \|(T - T_n)x_n^{(j)}\| + \|T_n\|\|x_n^{(i)} - x_n^{(j)}\| \qquad (3.5)$$

and, one gets by taking subsequences $\{z_i\} \subset \{x_n^{(i)}\}$, $\{z_j\} \subset \{x_n^{(j)}\}$

$$\|T z_i - T z_j\| \leq \|T - T_n\| (\|z_i\| + \|z_j\|) + \|T_n z_i - T_n z_j\| \qquad (3.6)$$

Since $\{T_n\} \to T$, we can find $n_0, i_0 \in N_0$ such that for $n(\geq n_0) \in N_0$, $\min(i,j) > i_0$, we have:

$$\|T - T_n\| < \frac{\varepsilon}{4c}, \quad \max(\|z_i\|, \|z_j\|) \leq K_z \leq 2c < \infty \text{ and } \|T_n z_i - T_n z_j\| < \frac{\varepsilon}{2}$$

for any given $c$ and $\varepsilon = \varepsilon(c) \in R_+$, since $\{z_i\}$ and $\{z_j\}$ are bounded subsequences, and $\{T_n z_i\}$ converges, so that it is a Cauchy sequence, since $\{T_n\}$ contains at least one compact operator. As a result, $\|T z_i - T z_j\| \leq \frac{\varepsilon}{2} + \frac{\varepsilon}{2} = \varepsilon$ is arbitrarily small for $\varepsilon$ being sufficiently small. Thus, $\{T z_i\}$ is convergent. Property (i) has been proven. Property (ii) follows from Property (i) and the fact that any operator composite sequence of bounded operators is a compact operator if there is at least one which is compact. □

Now, define the composite operator $\hat{T}(k+i+1, k): X \to X$; $\forall i, k(\geq i) \in N_0$ by

$$\hat{T}(k+i+1, k) = T_{k+i} \ldots T_{k+1} T_k = \left(T_{j_{k+i} k+i} \ldots T_{j_{k+i}-1, k+i} T_{1, k+i}\right) \ldots \left(T_{j_k k} \ldots T_{2k} T_{1k}\right); \forall x \in Dom(T_{10}) \quad (3.7)$$

; $\forall k \in N_0$. Define also the sequence $\{\hat{T}^0(k+i+1, k)\}$ of composite operators as $\hat{T}^0(k+i+1, k) = T^0_{k+i} \ldots T^0_{k+1} T^0_k$; $\forall k \in N_0$ where $T^0_{k+i}$ replaces each operator in the composite operator $T_{k+i}$ by its limit when such a limit exists. A result is now given based on the existence of the following limit:

$$\lim_{k \to \infty} d\left(\hat{T}(k+i+1, i)x, \hat{T}^0(k+i+1, i)x\right) = 0; \forall x \in Dom(T_{10})$$

The following result can be proven by using Lemmas 3.1-3.4:

**Theorem 3.5.** Consider the operator composite sequence $\{\hat{T}(k+i+1, k)\}$; $\forall k, i \in N_0$ of composed linear bounded operators in (2.15) on a Banach space $(X, \|\ \|)$, subject to $Im(T_{jk})(\neq \varnothing) \subset Dom(T_{j+1, k})$ for $j \in \overline{j_k}$, $Im T_{j_k k}(\neq \varnothing) \subset Dom T_{1, k+1}$ and $1 \leq j_k \leq J(\in N_0) < \infty$; $\forall k \in N_0$, and the sequence of composed linear operators $\{\hat{T}^0(k+i+1, k)\}$ of (2.14) defined in the same way as $\{\hat{T}(k+i+1, k)\}$ as $k \to \infty$ by replacing each operator possessing a limit by such a limit. The following properties hold:

**(i)** Either the sequences $\{\hat{T}(k+i+1, k)\}$ and $\{\hat{T}^0(k+i+1, k)\}$ have limits and both limits coincide or none of them has a limit and, furthermore, and $\hat{T}(k+i+1, k) \to \hat{T}^0(k+i+1, k)$ as $k \to \infty$.

**(ii)** If the limits of Property (i) exist and are finite then the limits of the sequences of operators $\{\hat{T}(k+i+1, k)\}$ and $\{\hat{T}^0(k+i+1, k)\}$ as $k \to \infty$; $\forall i \in N_0$ have the same set of fixed points.

**(iii)** Assume, in addition, that for some $k \in N_0$, there is at least one compact operator in the composition operator $\hat{T}(k+i+1,k)$, and that $Im(T_{jk})(\neq \emptyset) \subset \overline{Im}(T_{jk}) \subset Dom(T_{j+1,k})$ for $j \in \overline{j_k}$, $ImT_{j_k+\ell k}(\neq \emptyset) \subset \overline{Im}T_{j_k+\ell k}(\neq \emptyset) \subset DomT_{1,k+1}$ for $0 \leq \ell \leq i$ and some $i \in N_0$ and that all the operators are closed. If Property (i) holds with $\hat{T}(k+i+1,k) \to \hat{T}^0(k+i+1,k) \to \hat{T}^*$ as $k \to \infty$; for some $i \in N_0$ and $\|\hat{T}^*\| \leq K < 1$, then $\hat{T}^*$ is contractive and

$$\lim_{k,n\to\infty} \|\hat{T}(k+n(i+1),k)x - \hat{T}((k+n(i+1),k),k)y\| = \lim_{k\to\infty} \|\hat{T}^0(k+n(i+1),k)x - \hat{T}^0(k+n(i+1),k)y\|$$

$$= \lim_{k,n\to\infty} \|\hat{T}^{*(i+1)n}x - \hat{T}^{*(i+1)n}y\| = 0; \quad \forall x,y \in Dom(T_{10}) \qquad (3.8)$$

, the sequences of composite operators $\{\hat{T}(k+j,j)\}$ and $\{\hat{T}^0(k+j,j)\}$ converge to zero as $k \to \infty$, and $\hat{T}^*: Dom(\hat{T}^*) \subset X \to Im(\hat{T}^*) \subset X$ is bounded, closed and compact and has a unique fixed point in $Dom(\hat{T}^*) \cap \overline{Im}(\hat{T}^*)$ to which all sequences with initial conditions in $\overline{Im}(\hat{T}^*)$ converge.

**(iv)** Assume that there is a (in general, non unique) strictly increasing sequence of nonnegative integers $\{j_k\}$ with $j_0 = 0$ and $0 < j_{k+1} - j_k \leq m < \infty$ such that

$$\|\hat{T}(j_{k+2}, j_k)\| \leq \hat{K}(j_{k+2}, j_{k+1}) \|\hat{T}(j_{k+1}, j_k)\| \leq K \|\hat{T}(j_{k+1}, j_k)\| \qquad (3.9)$$

; $\forall k \in N_0$ for some nonnegative real sequence $\{K_k(j_{k+1}, j_k)\}; k \in N_0$ and some real constant $K \in [0,1)$. Assume, in addition, that for some $k \in N_0$, there is at least one compact operator in any composite operator $\hat{T}(j_{k+1}, j_k)$ and that all the operators are closed. Then, the sequences of composite operators $\{\hat{T}(k+j,j)\}$ and $\{\hat{T}^0(k+j,j)\}$ converge to zero as $k \to \infty$ for any finite $j \in N_0$. Finally, assume that $\hat{T}(j_{k+1}, j_k) \to \hat{T}_g^*$ as $k \to \infty$. Then, $\hat{T}_g^*$ is contractive, continuous, bounded, closed and compact and has a unique fixed point in $Dom(\hat{T}^*) \cap \overline{Im}(\hat{T}^*)$ to which all sequences with initial conditions in $\overline{Im}(\hat{T}^*)$ converge.